\documentclass[12pt,a4paper]{amsart}
\usepackage{amssymb}

\newtheorem{theo}{Theorem}
\newtheorem{lemma}[theo]{Lemma}
\newtheorem{proc}[theo]{Procedure}

\newtheorem{defi}[theo]{Definition}

\theoremstyle{definition}

\def\deg{\operatorname{deg}\nolimits}

\def\CC{{\mathcal{C}}}

\def\CH{{\mathcal{H}}}

\def\BC{{\mathbb{C}}}

\def\BP{{\mathbb{P}}}

\def\BQ{{\mathbb{Q}}}
\def\BZ{{\mathbb{Z}}}

\def\Aut{\operatorname{Aut}\nolimits}

\def\ie{{\em i.e.}}

\title{Variations on Van Kampen's method}
\author{David Bessis}
\address{David Bessis, LIFR-MI2P, CNRS / Independent University of Moscow,
Bolshoi Vlasevsky Pereulok, Dom 11, Moscow 121002, Russia.}
\email{bessis@mccme.ru}

\begin{document}

\begin{abstract}
We give a detailed account of the classical Van Kampen method for 
computing presentations of fundamental groups of complements
of complex algebraic curves, and of a variant of this method,
working with arbitrary projections (even with vertical asymptotes).
\end{abstract}

\maketitle    

{\bf Introduction}
In the 1930's, Van Kampen described a general technique for computing
presentations of fundamental groups of complements of complex
algebraic curves. Though Van Kampen's original approach was
essentially valid, some technical details were not entirely
clear and were later reformulated in more modern and rigorous terms
(see for example the account by Ch\'eniot, \cite{cheniot}).

It is possible to transform Van Kampen's ``method'' into an entirely
constructive algorithm. To my knowledge, two implementations have
been realized, one by Jorge Carmona, the other by Jean Michel and
myself (GAP package VKCURVE, \cite{VKCURVE}).

The goal of the present note is to clarify some aspects which are
usually neglected but must
be addressed to obtain an efficient implementation. Also, the ``Van
Kampen's method'' explained here differs from the classical one, which
assumes the choice of a ``generic'' projection: our variant method
works with an arbitrary projection. The reason for what may
appear to be a superfluous sophistication (since ``generic'' projections
always exists and are easy to find) is that working with a non-generic
projection may be computationally more efficient.
The variant method explained here is implemented in VKCURVE,
and has already been used to find previously unknown presentations.

\smallskip

Let
$P\in\BC[X,Y]$. The equation $P(X,Y)=0$ defines an algebraic
curve $\CC\in\BC^2$. Our goal is to find a presentation for the
fundamental group of $\BC^2-\CC$ (the method can be 
adapted to work with projective curves,
as it is briefly mentioned at the end of section
\ref{secmain}). Without loss of generality, we may (and will) assume
that $P$ is quadratfrei. View $P$ as a polynomial in $X$ depending
on the parameter $Y$:
$$P = \alpha_0(Y) X^d + \alpha_1(Y) X^{d-1} + \dots + \alpha_d(Y),$$
with $\alpha_0 (Y)\neq 0$.
To study $\BC^2-\CC$, we decompose it according to the fibers of the
projection $\BC^2 \rightarrow \BC,
(x,y)\mapsto y$. Up to changing the variables, one could
assume that $d$ equals the total degree of $P$ (the projection
 is then said to be ``generic''); however, for reasons detailed
below, we do not make this assumption that we have a generic projection.
For all but a finite number of exceptional values for $y_0$, the equation
$P(X,y_0)=0 $ has exactly $d$ distinct solutions in $X$.
The main idea in Van Kampen's method is that to compute a presentation,
it essentially suffices to be able to track these $d$ solutions when
the parameter $Y$ varies along certain loops (around the 
exceptional values). These $d$ solutions
form certain braids with $d$ strings, called {\em monodromy braids}.

In complexity terms, the most expensive part of the algorithm is
the computation of the monodromy braids. Since the computation
time increases with the number of strings, it is tempting to 
keep it reduced by working with non-generic projections, with a smaller number
of strings. However, these projections frequently involve vertical
asymptotes (to us, vertical lines are lines with equations of the form
$Y=y_0$, thus our $X$-axis is vertical, and our $Y$-axis is horizontal;
vertical asymptotes can only appear when $\alpha_0$ is not a scalar).
Classical Van Kampen method
is not adapted to deal with vertical asymptotes but,
as we explain in sections \ref{seclift} and \ref{monodromy}, some
corrections can be introduced to make it work.

The structure of this note is as follows: after some preliminaries in
section \ref{prelim}, we describe the main steps of the algorithm
in section \ref{secmain}. Sections \ref{seclift} and \ref{monodromy}
are devoted to two technical points -- they contain the only original
material of this note.

This note does not cover all aspects of the effective implementation
of Van Kampen's method. The most serious gap is that we do not
explain how to perform step $a$ of Procedure \ref{step3}; 
this, with other features of VKCURVE, will be described in a forthcoming
joint paper with Jean Michel.

{\bf \flushleft Note.}
A modification of Van Kampen's
method dealing with vertical asymptotes is also proposed in \cite{artal}.
Our approach is probably more or less equivalent to theirs, but formulated in
a way which is closer to a fully automated
procedure. The author thanks Jorge Carmona for useful electronic
discussions.

\section{Topological preliminaries}
\label{prelim}

Though we will use them only in dimension $2$, we formulate
the results of this section in arbitrary dimension, since it doesn't
cost more; for the same reason, we work in the projective space.

Let $n$ be a positive integer, let $\CH$ be an algebraic hypersurface
in the complex projective space $\BP^n$.

Let $x_0$ be a basepoint in $\BP^n - \CH$. We are interested in generating
$\pi_1(\BP^n-\CH,x_0)$. There is a natural class of elements
of this groups, the {\em meridians} (also called {\em
generators-of-the-monodromy}), in which to pick the generators.
To construct a meridian, one needs to choose:
\begin{itemize}
\item[1)] a smooth point $x\in \CH$;
\item[2)] a path from $x_0$ to $x$, intersecting $\CH$ only at the
endpoint $x$.
\end{itemize}
To these choices, one associates a loop as follows: start from $x_0$,
follow $\gamma$; just before reaching $x$, make a positive full
turn around $\CH$ (by the local inversion theorem,
in the neighbourhood of a smooth point, the complement
of an hypersurface ``looks like'' the complement of an hyperplane,
and the local fundamental group is isomorphic to $\BZ$ -- the standard
orientation of $\BC$ telling which generator is the positive one);
return to $x_0$ following $\gamma$ backwards.

The reader should convince himself that this makes sense and that
we have defined an element $s_\gamma \in \pi_1(\BP^n-\CH,x_0)$.
Of course, different choices may yield the same element.
However, although the element $x\in\CH_{\text{smooth}}$ is not
uniquely determined by $s_\gamma$, it should be noted that
it belongs to exactly one  of irreducible components
of $\CH$ (since intersections of components
belong to the singular locus), and that this component $D$
is uniquely determined by $s_\gamma$ (to see it, integrate over
$s_\gamma$ the inverses of defining polynomials).
We will say that $s_\gamma$ is a {\bf meridian of $\CH$ around $D$}.
It is important to note this notion depends not only on $D$, but
also on the remaining components of $\CH$, since these have to be
avoided when choosing the path from $x_0$ to a point of $D$.

We will get rid of all topological technicalities by
admitting without proof the following folk-lemma:

\begin{lemma}
\label{dust}
Let $D$ be an irreducible component of $\CH$. We have $\CH = \CH'\cup D$,
with $D\nsubseteq \CH'$, where $\CH'$ is the union of the remaining
components.
\begin{itemize}
\item[(i)] The meridians of $\CH$ around $D$ form a single conjugacy class.
\item[(ii)] Consider the embedding $\BP^n-\CH \hookrightarrow \BP^n-\CH'$,
and the associated morphism $\phi$ between fundamental groups.
Then $\phi$ is surjective, and its kernel is generated by
the meridians of $\CH$ around $D$.
\end{itemize}
\end{lemma}

One may already observe that this lemma has a meaning in terms of
generators and relations. First, an induction from (ii) proves that
$\pi_1(\BP^n-\CH)$ is generated by meridians. Secondly, assume that
we already know a presentation of $\pi_1(\BP^n-\CH)$, with generators
corresponding to meridians; then, by simply forgetting (as generators,
and in the relations) those generators which are meridians around
$D$, we obtain a presentation of  $\pi_1(\BP^n-\CH')$. We will use
this later.

One may also note that the complement of $m$ points in $\BC$ 
is relevant to the above discussion (with $n=1$).
Since we will need it the next sections, let us
fix the following {\em ad hoc} terminology:

\begin{defi}
\label{planartree}
Let $x_1,\dots,x_m$ be $m$ distinct points in $\BC$.
Let $x_0\in \BC-\{x_1,\dots,x_m\}$.
We define a
{\em planar tree connecting $x_0$ to $\{x_1,\dots,x_m\}$} to be a subset of
$\BC$ homeomorphic to a tree, containing $\{x_0,\dots,x_m\}$ and
such that $x_1,\dots,x_m$ are leaves.

Assume we have fixed such a planar tree $T$. For each $i\in\{1,\dots,m\}$,
there is (up to reparametrization) a unique path in $T$ connecting $x_0$
to $x_i$; it avoids the other leaves $(x_j)_{j\neq i}$. Therefore it
defines a meridian $s_i\in\pi_1(\BC-\{x_1,\dots,x_m\},x_0)$.
The {\em generating set associated to $T$} is, by definition,
$\{s_1,\dots,s_m\}$.
\end{defi}

It is clear, in the above definition, that the $s_i$'s
realize an explicit isomorphism between
$\pi_1(\BC-\{x_1,\dots,x_m\},x_0)$ and the free group on $m$
generators.

\section{The main idea : a fibration argument}
\label{secmain}

After these preliminaries, we move to the central matter.
Let $$P = \alpha_0(Y) X^d + \alpha_1(Y) X^{d-1} + \dots + \alpha_d(Y)$$
be as in the introduction.
For a ``generic'' choice of $y_0$, the equation $P(X,y_0)=0$ has $d$
solutions in $X$. Here ``generic'' means that $y_0$ should not
be a zero of the discriminant $\Delta \in \BC[Y]$ of $P$.

Let $y_1,\dots,y_r$ be the distinct roots of $\Delta$.
Let $B:=\BC-\{y_1,\dots,y_r\}$. Let
$E:=\{(x,y)\in \BC\times B | P(x,y) \neq 0\}$.
Using the classical fact that the roots of a polynomial
are continuous functions
of its coefficients, we see that the map $p: E\rightarrow B, (x,y) \mapsto y$
is a locally trivial fibration, with fibers homeomorphic to the complement
of $d$ points in $\BC$.

Choose a basepoint $y_0 \in B$. Let $F$ be the fiber over $y_0$.
In $F$, choose a basepoint $x_0$ (this choice is not innocent --
we will return to this in the next section).
Since $F$ is connected and $\pi_2(B)=1$, the fibration exact sequence
basically amounts to 
$$1 \longrightarrow \pi_1(F,x_0) \longrightarrow \pi_1(E,(x_0,y_0)) 
\longrightarrow \pi_1(B,y_0) \longrightarrow 1.$$

Both $\pi_1(F,x_0)$ and $\pi_1(B,y_0)$ are isomorphic to free groups,
of ranks respectively $d$ and $r$. Any short exact sequence landing
on a free group is split, so $\pi_1(E,(x_0,y_0))$ must be a semi-direct
product $F_d \rtimes F_r$.

Let us be more specific. Choose a planar tree connecting
$x_0$ to the roots of $P(X,y_0)$.
This spider defines free generators $f_1,\dots,f_d$
of $\pi_1(F,x_0)$, each of them being a meridian around one of the roots.

Similarly, choose a planar tree connecting
$y_0$ to the roots of $\Delta$. We obtain $r$ meridians $g_1,\dots,g_r$
generating freely $\pi_1(B,y_0)$.

Whichever way we choose to lift $g_1,\dots,g_r$ to elements $\tilde{g}_1,
\dots,\tilde{g}_r\in \pi_1(E,(x_0,y_0))$, we have a semi-direct
product structure. The conjugacy action
$f_i \mapsto \tilde{g}_j^{-1} f_i \tilde{g}_j$ defines an
{\em monodromy automorphism} $\phi_j \in \Aut (F_d)$ (we will explain in section
\ref{monodromy} how
to compute these automorphisms).
We obviously have the presentation
$$\pi_1(E,(x_0,y_0)) \simeq \left< f_1,\dots,f_d,\tilde{g}_1,\dots,
\tilde{g}_r | \tilde{g}_j^{-1} f_i \tilde{g}_j = \phi_j (f_i) \right>$$
(where by $\phi_j (f_i)$ we mean the corresponding word in the
$f_1,\dots,f_d$ and their inverses).

For each root $y_j$ of $\Delta$, define $L_j$ to be the line in $\BC^2$
of equation $Y=y_j$.
Clearly,
$$E = \BC^2 - \CC\cup L_1 \cup \dots \cup L_r.$$
Nothing prevents some of the $L_j$'s to be included in $\CC$.
The space we are interested in, $\BC^2-\CC$, is obtained from $E$
by adding (or, more exactly, by forgetting to remove) the $L_j$'s which
are not included in $\CC$.
What we are tempted to do is to use Lemma \ref{dust} (ii), and
to forget the corresponding $\tilde{g}_j$'s in the above presentation.
If we had chosen the $\tilde{g}_j$'s to be meridians, we would obtain
a presentation for $\pi_1(\BC^2-\CC)$.
The following Lemma proves that this strategy works:

\begin{lemma}
\label{lift}
Let $g_j\in \pi_1(B,y_0)$ be a meridian around $y_j$.
There exists in $\pi_1(E,(x_0,y_0))$ a meridian $\tilde{g}_j$
around $L_j$ such $p_*(\tilde{g}_j) = g_j$.
\end{lemma}

\begin{proof} This is a particular case of \cite{zardiag}, Lemma 2.4.
(whose proof relies on an easy general position argument).
For a more constructive argument, see next section.
\end{proof}

Van Kampen's method is the following procedure, which summarizes the above
discussion:

\begin{proc}[``Van Kampen's method'']
\label{main}
Start with a (quadrafrei) polynomial $P\in \BC[X,Y]$.
\begin{itemize}
\item[1)] Compute $\Delta \in \BC[Y]$ and its roots $y_1,\dots,y_r$.
\item[2)] Choose a basepoint $y_0\in \BC-\{y_1,\dots,y_r\}$.
Construct generating
meridians $g_1,\dots,g_r$ in $\pi_1(\BC-\{y_1,\dots,y_r\},y_0)$.
Choose $x_0\in \BC$ such that $P(x_0,y_0)\neq 0$. Lift
$g_1,\dots,g_r$ to elements $\tilde{g}_1,\dots,\tilde{g}_r\in
\pi_1(E,(x_0,y_0))$ which are meridians around the $L_j$'s.
\item[3)] Choose generating meridians $f_1,\dots,f_d$ for the fundamental
group of the fiber over $y_0$.
Compute the monodromy automorphisms $\phi_1,\dots,\phi_r$.
\item[4)] Let $J:=\{j\in\{1,\dots,r\} | L_j \nsubseteq \CC \}$. 
We obtain the following presentation for $\pi_1(\BC^2-\CC,(x_0,y_0))$:
\small{
$$\left< f_1,\dots, f_d, \tilde{g}_1,\dots,\tilde{g}_r \left| 
\begin{array}{c}
\forall i \in \{1,\dots,d\}, \forall j \in \{1,\dots,r\},

\tilde{g}_j^{-1} f_i \tilde{g}_j = \phi_j (f_i) \\
\forall j\in J, \tilde{g}_j = 1
\end{array}
\right.
\right>$$}
\end{itemize}
\end{proc}

In implementations of step $1$, one only computes approximations
of $y_1,\dots,y_r$; it is not difficult to make sure they are
good enough for our purposes (which is to find loops circling
the actual roots).
In the next two sections, we will describe explicit constructive
versions of steps $2$ and $3$.
To be really useful, step $4$ should be followed by a procedure
``simplifying'' the initial presentation, which is usually highly
redundant and unpleasant. This problem has no general solution,
though some heuristics, implemented in VKCURVE, happen to be
quite effective on many practical examples; these heuristics will
not be described here.

{\flushleft \bf Projective Van Kampen method}
The above procedure explains how to compute the fundamental
group of the complement of an affine complex algebraic curve.
To deal with a projective curve $\CC\subset \BC\BP^2$
proceed as follows: decompose $\BC\BP^2$ as $\BC^2\cup \BC\BP^1$; let
$\CC':=\CC\cap\BC^2$ with
the above method, compute a presentation for $\BC^2-\CC' =
\BC\BP^2 - \CC\cup \BC\BP^1$; let $s_\infty\in\pi_1(\BC\BP^2-\CC\cup\BP^1)$
be a meridian of $\CC\cup \BC\BP^1$ around $\BP^1$ (in most situations,
if $f_1,\dots,f_d$ are as in section \ref{monodromy},
$s_\infty:=(f_df_{d-1}\dots f_1)^{-1}$ is suitable;
see \cite{artal}, Prop. 2.10); by Lemma \ref{dust},
by adding the relation $s_\infty=1$, one obtains a presentation
for $\pi_1(\BC\BP^2-\CC)$.

\section{Lifting meridians}
\label{seclift}

In this section, we discuss the problem of lifting generating meridians
$g_1,\dots,g_r$ of $\pi_1(B,y_0)$ to elements $\tilde{g}_1,\dots,
\tilde{g}_r$ of $\pi_1(E,(x_0,y_0))$. As explained in the previous
section, we would like these $\tilde{g_j}$'s to be meridians around the
$L_j$'s. Thanks to Lemma \ref{lift}, we know that this is possible. However,
we would like to do this in a constructive manner -- something we could
easily instruct a computer to do.

Practically, to encode the $g_j$'s, we choose representants
$\gamma_1,\dots,\gamma_r$ in the loop space $\Omega(B,y_0)$
(in VKCURVE, we actually work with piecewise-linear
loops, with endpoints in $\BQ[i]$).

Thus the most natural way of lifting the $g_j$'s is to take the
elements represented by the loops $\tilde{\gamma}_1,\dots,\tilde{\gamma}_r$,
defined as follows: for all $t\in [0,1]$ and all $j$, we set
$\tilde{\gamma}_j(t):=(x_0,\gamma_j(t))$. Two possible problems arise from
this idea:
\begin{itemize}
\item
The first problem is that $\tilde{\gamma}_j$ may not be a loop
in $E$: nothing prevents it from intersecting $\CC$ (although, by
general position arguments, it should only happen for a finite
number of unlucky choices for $x_0$).
\item
The second problem is a more serious one: the $\tilde{\gamma}_j$'s
may not represent meridians.
\end{itemize}

The classical Van Kampen method assumes that $\CC$ does not have 
vertical asymptotes. With this assumption, a compactness argument
can be used to eliminate both problems: the supports of the $\gamma_j$'s
can be assumed to be all included in a disk $D$, and, if we choose
$x_0$ large enough (in module), we can guarantee that $\CC$ does
not intersect the lifted disk $(x_0,D)$. With such a choice,
it follows easily that the $\tilde{\gamma}_j$'s indeed represent meridians
around the $L_j$'s.

However, as explained in the introduction, we do not want to assume
the absence of vertical asymptotes. To ensure that none of the above two
problems occurs, we rely on the following criterion
(as in Procedure \ref{main}, we set $J:=\{j\in\{1,\dots,r\} | L_j
\nsubseteq \CC\}$):

\begin{lemma}[Explicit Step 2]
\label{step2}
Let $Q:=P\prod_{j\in J} (Y-y_j)$.
Let $\nabla\in\BC[X]$ be the discriminant of $Q$, viewed as polynomial
in $Y$ with coefficients in $\BC[X]$.
Choose $x_0\in \BC$ which is not a root of $\nabla$.
Denote by $S$ the set of solutions in $Y$ of $Q(x_0,Y)=0$.
For all $j\in \{1,\dots,r\}$, one has $y_j\in S$.

Choose $y_0\in \BC - S$, and for all $j\in \{1,\dots,r\}$ choose a loop
$\gamma_j \in \Omega(\BC-S,y_0)$ representing a meridian of $S$ around
$y_j$. Then, for all $j$, the path
$\tilde{\gamma}_j:= (x_0,\gamma_j)$ represents
a meridian of $\CC\cup L_1 \cup \dots \cup L_r$ around $L_j$.
\end{lemma}

\begin{proof}
For all $j\in\{1,\dots,r\}$, either $j\in J$, or $L_j\subseteq \CC$;
in both case, it is clear that $y_j \in S$.
The path $\tilde{\gamma}_j$ avoids $\CC\cup L_1\cup \dots \cup L_r$,
since the intersection of this curve with the line $X=x_0$ is 
precisely described by $S$.
Proving that $\tilde{\gamma}_j$ represents a meridian of
$\CC\cup L_1\cup \dots \cup L_r$ around $L_j$
essentially amounts to checking that
$(x_0,y_j)$ is a smooth point of $\CC\cup L_1\cup \dots \cup L_r$ or,
by an immediate reformulation, that $L_j$ is the only component
of $\CC\cup L_1\cup \dots \cup L_r$ in which $(x_0,y_j)$ lies.
If this was not satisfied, the polynomial $Q(x_0,Y)$
would have multiple roots in $Y$, which contradicts the assumption on $x_0$.
\end{proof}

\section{Computing the monodromy automorphisms}
\label{monodromy}

At the end of the second step of Procedure \ref{main}, as detailed
in Lemma \ref{step2}, we are provided with:
\begin{itemize}
\item A basepoint $(x_0,y_0)\in E$.
\item Loops $\gamma_1,\dots,\gamma_r \in \Omega(B,y_0)$ such
that each horizontally lifted loop $\tilde{\gamma}_j=(x_0,\gamma_j)$
is in $\Omega(E,(x_0,y_0))$ and represents a meridian
of $\CC\cup L_1\cup \dots \cup L_r$ around $L_j$.
\end{itemize}

As it was explained in the previous section,
in the classical Van Kampen method where one assumes the absence of vertical
asymptotes, one can get rid of many problems by choosing $x_0$ far enough;
then, to compute the {\em monodromy automorphism} $\phi_j$ corresponding to
$\tilde{\gamma}_j$, it suffices
to track the solutions in $X$ of $P(X,Y)$ when $Y$ moves along $\gamma_j$:
this defines a {\em monodromy braid} $b_j$ on $d$ strings, from which
the automorphism $\phi_j$ can be deduced using the standard Hurwitz
action of the braid group on the free group
(see Definition \ref{defihurwitz} and Lemma \ref{hurwitz}).

However, since we decided
to work in a situation allowing vertical asymptotes, we may no longer
assume that $x_0$ is ``far enough''; in particular, it may
occur that the strings of the monodromy braid turn around $x_0$, in
which case the monodromy automorphism cannot be computed by
Hurwitz formulas (the example given at the end of this section should
convince the reader that there is a serious obstruction -- this is not
just a matter of being smart when choosing $x_0$).
However, by adding to the monodromy braid an additional
string
fixed at $x_0$, one obtain extra information which can be used
to modify Hurwitz formulas in a suitable way.
This is what we detail in the present section.

To simplify notation, we fix some $j\in \{1,\dots,r\}$, and write
$\gamma$ and $\tilde{\gamma}$ instead of $\gamma_j$ and $\tilde{\gamma}_j$;
our goal is to compute the corresponding monodromy automorphism $\phi$.

Let $\{x_1,\dots,x_{d+1}\}$ be the set of solutions in $X$ of
$(X-x_0)P(X,y_0)=0$. For the sake of simplicity, we
assume that the $x_i$'s have distinct real parts (this is always true
up to rescaling), and that $\Re(x_1) < \Re(x_2) < \dots < \Re(x_{d+1})$
(this is always true up to reordering). Among the $x_i$'s is $x_0$, say
$x_0  = x_{i_0}$.

Let $X_{d+1}$ be the configuration space of $d+1$ points in the complex
line. Tracking the solutions of $(X-x_0)P(X,\gamma(t))=0$ for $t\in [0,1]$,
we obtain an element of $\Omega(X_{d+1},\{x_1,\dots,x_{d+1}\})$, which
represents an element $b_\gamma$ of the braid group $\pi_1(X_{d+1},
\{x_1,\dots,x_{d+1}\})$.

For all $i\in \{1,\dots,d+1\}$, we denote by $x_i(t)$ the string
of the monodromy braid starting at $x_i$, \ie, the
unique continuous path $[0,1]\rightarrow \BC$ such that
$\forall t\in [0,1], P(x_i(t),\gamma(t))=0$ and $x_i(0)=x_i$.

Let $F:=\BC-\{x_1,\dots,x_{i_0-1},x_{i_0+1},\dots,
x_{d+1}\}$. We view $x_0$ as a basepoint for $F$.
We also introduce a secondary basepoint $x_\infty$, which will be
used in our argumentation but will not appear in the formulation of the
final result. This secondary basepoint is assumed to have a ``negative
enough'' imaginary part,
in the sense that it satisfies the following conditions:
$$\forall t\in [0,1], \forall i \in \{1,\dots, d+1\}, \Im(x_\infty) <
\Im(x_i(t))$$
and
$$\forall i \in \{1,\dots, d\}, \Re
\left(\frac{x_{i+1}-x_\infty}{x_i-x_\infty}\right) >0.$$
It is not difficult to figure out why such $x_\infty$ exist; let us fix one.

We consider the planar tree $\bigcup_{i=1}^{d+1}[x_\infty,x_i]$
(where $[x_\infty, x_i]$ denotes the linear segment between $x_\infty$
and $x_i$). We use this tree to describe generators for
various fundamental groups (see Definition \ref{planartree}):
\begin{itemize}
\item Being a tree connecting $x_\infty$ to
 $\{x_1,\dots,x_{d+1}\}$, it defines generating meridians
 $$e_1,\dots,e_{d+1}$$ for $\pi_1(F-\{x_0\},x_\infty)$.
\item 
Being a tree
connecting $x_\infty$ to
$\{x_1,\dots,x_{i_0-1},x_{i_0+1},\dots,
x_{d+1}\}$, it defines generating meridians
for $\pi_1(F,x_\infty)$. Conveniently abusing notations, we still
denote them by $$e_1,\dots,e_{i_0-1},e_{i_0+1},\dots,e_{d+1}.$$
\item Being a tree connecting $x_0=x_{i_0}$ to
$\{x_1,\dots,x_{i_0-1},x_{i_0+1},\dots,
x_{d+1}\}$, it defines generating meridians
$$f_1,\dots,f_{i_0-1},f_{i_0+1},\dots,f_{d+1}$$

for $\pi_1(F,x_0)$.
\end{itemize}

Let $E_\gamma$ be the pullback over
$[0,1] \stackrel{\gamma}{\rightarrow} B$ of the fiber bundle
$E \stackrel{p}{\twoheadrightarrow} B$; in other words, the fiber $E_t$ of
$E_\gamma$ over $t$ is the complement in $\BC$ of $\{x_1(t),\dots,x_{i_0-1}(t),
x_{i_0+1}(t),\dots,x_{d+1}(t)\}$.
The space $F$ defined above coincides with the fiber over $0$ (or,
equivalently, $1$).

\begin{lemma}
\label{triv}
There exists a trivialization 
$$\Psi:E_\gamma \simeq F\times [0,1]$$ such that, for
all $t\in[0,1]$, $\Psi((x_0,\gamma(t)))= (x_0,t)$ and
$\Psi((x_\infty,\gamma(t)))= (x_\infty,t)$.
In particular, the induced homeomorphism
$$\psi:F=E_0\stackrel{\sim}{\rightarrow} E_1=F$$ satifies
$\psi(x_0)=x_0$ and $\psi(x_\infty)=\psi(x_\infty)$.
\end{lemma}

\begin{proof}
This is an elementary
variation of the standard construction of the map from the braid
group to the mapping class group of the punctured plane.
Basically, one has to imagine that the plane is a piece of rubber, pinned
to a desk at $x_0$ and $x_\infty$, and that we force $d$ other points
to move according to $b_\gamma$. We leave the details to the reader.
\end{proof}

We choose $\Psi$ and $\psi$ as in the lemma.
Since $\psi$ fixes both $x_0$ and $x_\infty$, it induces  automorphisms
$\psi_{*0}\in \Aut (\pi_1(F,x_0))$ and $\psi_{*\infty}\in
\Aut(\pi_1(F,x_\infty))$. 
We will need a third automorphism:
since $\psi(x_0)=x_0$, $\psi$ restricts to an homeomorphism
$\tilde{\psi}$ of $F-\{x_0\}$, and induces an element
$\tilde{\psi}_{*\infty} \in \Aut(\pi_1(F-\{x_0\},x_\infty))$.

Let $\tilde{\gamma}$ be the lifted path $(x_0,\gamma)\in\pi_1(E,(x_0,y_0))$.
Let $\phi\in\Aut(\pi(F,x_0))$ be the associated monodromy automorphism
(see section \ref{secmain}).

\begin{lemma}
We have $\phi=\psi_{*0}$.
\end{lemma}

\begin{proof}
Consider the loop $\tilde{\gamma}=(x_0,\gamma)\in\Omega(E,(x_0,y_0))$.
In the pull-back $E_\gamma$, it corresponds to the horizontal
path $t\mapsto (x_0,t)$.
For any loop $\omega\in\Omega(F,x_0)$, we may use the trivialization
of Lemma \ref{triv} to construct a homotopy in $\Omega(E_\gamma,(x_0,0))$
between
$\omega$ (viewed as a loop in the fiber of $E_\gamma$ over $0$)
and $\tilde{\gamma} \psi(\omega) \tilde{\gamma}^{-1}$ (where $\psi(\omega)$
is viewed as a loop in the fiber of $E_\gamma$ over $1$). Pushed
back in $E$, this homotopy shows that conjugating by $\tilde{\gamma}$ is
the same as applying $\psi_{*0}$.
\end{proof}

We may now explain our strategy for computing $\phi$.
First, we compute $\tilde{\psi}_{*\infty}$.
As announced earlier, since $x_\infty$
if ``far enough'', this can be done using Hurwitz formulas. Then
we use this intermediate result in two ways: first, we deduce
$\psi_{*\infty}$; then, we compute the discrepancy between
$\psi_{*\infty}$ and $\psi_{*0}$ coming 
from the change of basepoint.

Since this is the convenient setting for implementing the method,
we will suppose we are able to write $b_\gamma$ as a word in
the standard generators of the braid group:

\begin{defi}[Abstract braid group]
We denote by $B_{d+1}$ the group 
given by the abstract presentation
$$\left< \sigma_1,\dots,\sigma_{d} \left| 
\begin{array}{c}
\sigma_i\sigma_{i+1}\sigma_i 
=\sigma_{i+1}\sigma_i\sigma_{i+1} \text{ for all } i \\
\sigma_i\sigma_j = \sigma_j\sigma_i \text{ for all } i,j
\text{ with } i-j>1
\end{array} \right. \right>$$
\end{defi}
It is well-known that $B_{d+1}$ is isomorphic to the fundamental group
of $X_{d+1}$. Let us be more specific. For each $i\in\{1,\dots,d\}$,
the positive twist of two consecutive 
strings along the segment $[x_i,x_{i+1}]$
defines an element $\sigma_i\in\pi_1(X_{d+1},\{x_1,\dots,x_{d+1}\})$;
these elements satisfy the relations of the above definition and
realize an explicit isomorphism (this is nothing but the usual way
of considering braids via their real projection).

\begin{defi}[Hurwitz action]
\label{defihurwitz}
For all $i\in\{1,\dots,d\}$ and $j\in\{1,\dots,d+1\}$, set
$$H_i(e_j) := \left\{ \begin{array}{cc}
e_{i+1} & \text{ if } i=j \\
e_{i+1}e_ie_{i+1}^{-1}& \text{ if } i=j+1 \\
e_j & \text{ otherwise}
\end{array}
\right.$$
This defines an automorphism $H_i$ of the free group $\left<e_1,\dots,e_{d+1}
\right>$.
The $H_i$'s satisfy the defining relations of $B_{d+1}$ and induce
a morphism $H:B_{d+1} \rightarrow \Aut (\left<e_1,\dots,e_{d+1}\right>)$.
\end{defi}

{\bf \flushleft Note.} For compatibility with 
usual conventions for multiplication in
fundamental groups, we assume in the above definition and throughout this
section that groups of automorphisms act on the right.

\begin{lemma}
\label{hurwitz}
The automorphism $\tilde{\psi}_{*\infty}$ induced by
$\tilde{\psi}$ on
 $\left<e_1,\dots,e_{d+1}\right>
= \pi_1(F-\{x_0\},x_{\infty})$
 is $H(b_\gamma)$.
\end{lemma}

\begin{proof}
As in Lemma \ref{triv}, we can construct, for all 
braid $b\in B_{d+1}$, an homeomorphism of the 
pointed space $(F-\{x_0\},x_\infty)$. This
obviously induces a morphism $B_{d+1} \rightarrow \Aut(
\left<e_1,\dots,e_{d+1}\right> )$.
This morphism coincides with Hurwitz action
(it is enough to check this for the standard
generators of $B_{d+1}$, which is easy and classical).
The result follows as a particular case.
\end{proof}

Via the inclusion $F-\{x_0\} \hookrightarrow F$, we have
$$\pi_1(F,x_\infty) = \left< e_1,\dots,e_{d+1}\right> /e_{i_0}.$$ As mentioned
above, we will still denote by $e_i$ the image of $e_i$
in the quotient.
The braid $b_\gamma$ is $x_0$-pure (since the $x_0$-strand is constant);
therefore, the automorphism $H(b_\gamma)$ sends $e_{i_0}$ to a conjugate
of $e_{i_0}$. In particular, $H(b_\gamma)$ induces an endomorphism
of $\left< e_1,\dots,e_{d+1}\right> /e_{i_0}$, which is nothing but
the automorphism
of $\pi_1(F,x_\infty)$ induced by $\psi$.

An automorphism of a topological space yields a natural automorphism
of the functor from the fundamental groupoid to the category of groups,
which associates to each point the fundamental group at this point.
We have the following
commutative diagram of isomorphims, where the vertical arrows are
isomorphisms associated to paths connecting the two basepoints:

\[\begin{array}{cccc}
& \pi_1(F,x_0) & \xrightarrow{\phantom{mi}\psi_{*0}=\phi\phantom{mi}} & 
\pi_1(F,x_0) \\
\vcenter{\scriptsize{\rlap{$h_{[x_0,x_\infty]}$}}} &
\Big\downarrow & &
\Big\downarrow\vcenter{\scriptsize{\rlap{$h_{\psi([x_0,x_\infty])}$}}} \\
& \pi_1(F,x_\infty) & \xrightarrow{\phantom{mm}\psi_{*\infty}\phantom{mm}}
& \pi_1(F,x_\infty)
\end{array}\]

\smallskip

Our goal is to compute the monodromy automorphism
$$\phi = h_{\psi([x_0,x_\infty])}^{-1}\psi_{*\infty} 
h_{[x_0,x_\infty]} = 
(h_{\psi([x_0,x_\infty])}^{-1}  h_{[x_0,x_\infty]} )
(h_{[x_0,x_\infty]}^{-1}\psi_{*\infty}
h_{[x_0,x_\infty]}).$$
Since $h_{[x_0,x_\infty]}$ is the isomorphism
sending $f_i$ to $e_i$,
the automorphism $h_{[x_0,x_\infty]}^{-1}\psi_{*\infty}
h_{[x_0,x_\infty]}$ is given by Hurwitz action (Lemma \ref{hurwitz} --
of course after replacing the $e_i$'s by the $f_i$'s).

The automorphism $h_{\psi([x_0,x_\infty])}^{-1}  h_{[x_0,x_\infty]}$
is an inner automorphism. This is where we use the extra information
provided by the $x_0$-string of $b_\gamma$:

\begin{lemma}
\label{inner}
The element $H(b_\gamma)(e_{i_0})$ is conjugate to $e_{i_0}$ in
$\left< e_1,\dots,e_{d+1}\right>$. Let $a$ be such that
$H(b_\gamma)(e_{i_0}) = ae_{i_0}a^{-1}$. Let $\overline{a}$ be the
image of $a$ in $\left< f_1,\dots,f_{i_0-1},f_{i_0+1},\dots,f_{d+1}
\right>$ by the morphism sending $e_i,i\neq i_0$ to $f_i$ and $e_{i_0}$ to
$1$.
Then $h_{\psi([x_0,x_\infty])}^{-1}  h_{[x_0,x_\infty]}$ is the morphism
$f\mapsto \overline{a}^{-1}f\overline{a}$.
\end{lemma}

\begin{proof}
Left to the reader \small{(Hint: if true, this lemma provides a formula
for $\phi$; first, check that this formula indeed defines a morphism
$b_\gamma \mapsto \phi$; then check it on
generators of the group of $x_0$-pure braids on $d+1$ strings)}.
\end{proof}

The following procedure is a summary of the results of this section.
As promised, the exact choice of $x_\infty$ does not matter; 
nor does it matter to have distinct notations for the
$e_i$'s and the $f_i$'s.

\begin{proc}[Explicit Step 3]
\label{step3}
Suppose $\gamma$ is one the loops $\gamma_j$ constructed in 
Lemma \ref{step2}. To compute the associated monodromy automorphism
$\phi$, one may proceed as follows:
\begin{itemize}
\item[a)] Compute (as a word in the standard generators, using
real projection) the monodromy braid $b_\gamma$ with $d+1$ strings.
\item[b)] Compute Hurwitz action $H(b_\gamma)$ on the free group
$\left<f_1,\dots,f_{d+1}\right>$.
\item[c)] Identify the index $i_0$ of the $x_0$-string.
\item[d)] Find $a$ such that $H(b_\gamma)(f_{i_0})= af_{i_0}a^{-1}$
(this is trivial to do: take $a$ to be the first half of a reduced
word for $H(b_\gamma)(f_{i_0})$).
\item[e)] The composition of
$H(b_\gamma)$ with the automorphism $f\mapsto a^{-1}fa$ is
an automorphism of $\left<f_1,\dots,f_{d+1}\right>$ fixing $f_{i_0}$.
It induces an automorphism of the free group $\left<f_1,\dots,f_{d+1}\right>/
f_{i_0}$ of rank $d$: this is the monodromy automorphism $\phi$.
\end{itemize}
\end{proc}

Steps $b$, $c$, $d$ and $e$ are straightforward to implement, as soon as
one works with a software where braid groups, free groups and
groups automorphisms are available (this is the case with GAP).
Finding an efficient implementation of step $a$ is the main issue.

In the classical method, one may assume 
that $x_0$ has a ``large enough'' real part; this implies that $i_0=d+1$,
that
$H(b_\gamma)$ has no factor $\sigma_{d}^{\pm 1}$ and
that $H(b_\gamma)(f_{d+1})=f_{d+1}$:
Hurwitz action does not need a corrective term.

{\bf \flushleft Complexity of the modified method.} As mentioned in the
introduction, we are interested in non-generic projections because
they reduce the number of strings ($d < \deg P$).
To be able to work in this context, it is necessary
to introduce an additional string. The complexity cost from
this additional string is usually much smaller than the 
gain ($d+1\leq \deg P$; also, 
the additional string, which is constant at $x_0$, is
handled very efficiently). The cost of steps $c$, $d$ and $e$
(which are not present in classical Van Kampen method) is
neglectible. Actually, the only serious side-effect of the 
variant method is not in Procedure \ref{step3} itself,
but hidden in Lemma \ref{step2}: one has to construct more complicated
loops (for which step $3a$ will be more costly). However, the overall
balance is positive (this is a purely empirical statement,
we did not try to assess the theoretical complexity of our implementation).

\smallskip

To conclude, we illustrate by an example
(inspired by Example 3.1 in \cite{artal}) how the 
monodromy automorphism is affected by the choice of the basepoint.
It is a good exercise to apply Procedure \ref{step3} to check
the claims.

{\bf \flushleft Example.} Consider the example of the monodromy
braid with two strings whose positions are
$e^{-i\pi t}$ and $-e^{-i\pi t}$, for $t\in [0,1]$. This braid
occurs when studying the monodromy of $(XY-1)(XY+1)=0$ around the
singular fiber $Y=0$.
The induced monodromy automorphism depends on the choice of the
basepoint in $\BC-\{\pm1\}$:
\begin{itemize}
\item[i)] choosing $x_\infty:= -2i$ as basepoint for the fiber,
the monodromy automorphism has infinite order;
\item[ii)] choosing $x_0:=0$ as basepoint for the fiber, the monodromy
automorphism has order $2$. In particular, it cannot be described in
terms of usual Hurwitz action.
\end{itemize}

\end{document}